\documentclass{amsart}

\newtheorem{theorem}{Theorem}[section]
\newtheorem{lemma}[theorem]{Lemma}
\newtheorem{proposition}[theorem]{Proposition}

\theoremstyle{definition}

\newtheorem{example}[theorem]{Example}

\theoremstyle{remark}

\reversemarginpar

\begin{document}

\title[Powers, Exponentials and Logarithms]{The 
integrals in Gradshteyn and Ryzhik. \\
Part 8: Combinations of powers, exponentials and logarithms.}

\subjclass[2000]{Primary 33}


\author{Victor H. Moll}
\address{Department of Mathematics,
Tulane University, New Orleans, LA 70118}
\email{vhm@math.tulane.edu}

\author{Jason Rosenberg}
\address{Department of Mathematics,
Tulane University, New Orleans, LA 70118}
\email{jrosenbe@tulane.edu}

\author{Armin Straub}
\address{Department of Mathematics,
Tulane University, New Orleans, LA 70118}
\email{astraub@math.tulane.edu}

\author{Pat Whitworth}
\address{Department of Mathematics,
Tulane University, New Orleans, LA 70118}
\email{pwhitwor@tulane.edu}

\thanks{The authors  wish to the 
partial support of NSF-CCLI 0633223.}

\begin{abstract}
We describe some examples of integrals from the table of 
Gradshteyn and Ryzhik where the integrand is a combination of powers, 
exponentials and logarithms. The expressions for some of these integrals 
involve the Stirling numbers of the first kind.
\end{abstract}

\maketitle

\section{Introduction} \label{intro} 
\setcounter{equation}{0}

The uninitiated reader of the table of integrals by 
I. S. Gradshteyn and I. M. Ryzhik \cite{gr} will surely be 
puzzled by choice of integrands. In this note we provide an elementary proof 
of the evaluation {$\mathbf{4.353.3}$} 
\begin{equation}
\int_{0}^{1} ( ax + n+1)x^{n} e^{ax} \ln x \, dx = 
e^{a} \sum_{k=0}^{n} (-1)^{k-1} \frac{n!}{(n-k)! a^{k+1}} +
(-1)^{n} \frac{n!}{a^{n+1}}.
\label{iden1}
\end{equation}
\noindent
We also consider the integrals
\begin{equation}
q_{n} := \int_{0}^{1} x^{n} e^{-x} \ln x \, dx 
\end{equation}
\noindent
and the companion family 
\begin{equation}
p_{n} := \int_{0}^{1} x^{n} e^{-x} \, dx.
\end{equation}
\noindent
The integral $q_{n}$ corresponds to the case $a=-1$ in (\ref{iden1}). 
Section \ref{new} provides closed-form
expressions for $p_{n}$ and $q_{n}$. Section \ref{sec-param} considers the 
generalization
\begin{equation}
P_{n}(a) = \int_{0}^{1} x^{n} e^{-ax} \, dx \text{ and }
Q_{n}(a) = \int_{0}^{1} x^{n} e^{-ax} \ln x \, dx. 
\end{equation}
\noindent
The main result of this section is the closed-form expressions
\begin{equation}
P_{n}(a) := \int_{0}^{1}x^{n}e^{-ax} \, dx = 
\frac{n!}{a^{n+1}} \left( 1 - e^{-a} \sum_{k=0}^{n} \frac{a^{k}}{k!} \right),
\end{equation}
\noindent
and
\begin{equation}
Q_{n}(a) := \int_{0}^{1}x^{n}e^{-ax} \ln x \, dx = 
\frac{n!}{a^{n+1}} \left[  \sum_{k=1}^{n} \frac{1}{k} 
\left( 1 - e^{-a} \sum_{j=0}^{k-1} \frac{a^{j}}{j!} \right) + a Q_{0}(a) 
\right], 
\nonumber 
\end{equation}
where 
\begin{equation}
Q_{0}(a) = \int_{0}^{1} e^{-ax} \ln x \, dx = -\frac{1}{a} 
\left( \gamma + \ln a + \Gamma(0,a) \right),
\end{equation}
\noindent
and $\Gamma(0,a)$ is the incomplete gamma function defined by
\begin{equation}
\Gamma(a,x) := \int_{x}^{\infty} t^{a-1} e^{-t} \, dt.
\end{equation}

\section{The evaluation of {$\mathbf{4.353.3}$}} \label{simple} 
\setcounter{equation}{0}

The identity 
\begin{equation}
\frac{d}{dx} \left( x^{n+1}e^{ax} \right) = (ax+n+1)x^{n}e^{ax} 
\end{equation}
\noindent
and integration by parts yield 
\begin{equation}
\int_{0}^{1} ( ax + n+1)x^{n} e^{ax} \ln x \, dx = 
- \int_{0}^{1}x^{n} e^{ax} \, dx. \label{form2}
\end{equation}
\noindent
This last integral appears as $\mathbf{3.351.1}$ in \cite{gr}. We have obtained
a closed-form expression for it in \cite{moll-gr7}. A new proof is 
presented in Section \ref{sec-param}. 

A closed form expression for the right hand side of  (\ref{form2}) 
is obtained from 
\begin{equation}
\int_{0}^{1} x^{n} e^{ax} \, dx = \left( \frac{d}{da} \right)^{n} 
\frac{e^{a}-1}{a}. 
\label{form3}
\end{equation}
\noindent
The symbolic evaluation of (\ref{form3}) for small values of 
$n \in \mathbb{N}$ suggests the existence of a polynomial $P_{n}(a)$
such that 
\begin{equation}
\int_{0}^{1}x^{n}e^{ax} \, dx = 
\frac{(-1)^{n+1} \, n!}{a^{n+1}} + \frac{P_{n}(a)}{a^{n+1}} e^{a}. 
\label{pn-def}
\end{equation}
\noindent 
The next lemma confirms the existence of this polynomial. 

\begin{lemma}
The function $P_{n}(a)$ defined by 
\begin{equation}
P_{n}(a) =a^{n+1}e^{-a} \left( 
\left( \frac{d}{da} \right)^{n} \frac{e^{a}-1}{a} - 
\frac{(-1)^{n+1} n!}{a^{n+1}} \right) 
\end{equation}
\noindent
is a polynomial of degree $n$. 
\end{lemma}
\begin{proof}
Let $D = \frac{d}{da}$. Then $D^{n+1} = D(D^{n})$ produces the recurrence
\begin{equation}
P_{n+1}(a) = a P_{n}'(a) + (a-n-1)P_{n}(a).
\label{recu1}
\end{equation}
\noindent
The initial condition $P_{0}(a) = 1$ and (\ref{recu1}) show that 
$P_{n}$ is a polynomial of degree $n$. 
\end{proof}

\begin{theorem}
The polynomial 
\begin{equation}
Q_{n}(a) := (-1)^{n} P_{n}(-a)
\end{equation}
\noindent
has positive integer coefficients, written as 
\begin{equation}
Q_{n}(a) = \sum_{k=0}^{n} b_{n,k}a^{k}.
\end{equation}
\noindent
These coefficients satisfy
\begin{eqnarray}
b_{n+1,0} & = & (n+1)b_{n,0} \label{recu3} \\
b_{n+1,k} & = & (n+1-k)b_{n,k} + b_{n,k-1}, \quad 1 \leq k \leq n  \nonumber \\
b_{n+1,n+1} & = & b_{n,n}. \nonumber
\end{eqnarray}
\noindent
Moreover, the polynomial $Q_{n}(a)$ is given by
\begin{equation}
Q_{n}(a) = n! \, \sum_{k=0}^{n} \frac{a^{k}}{k!}
\end{equation}
\end{theorem}
\begin{proof}
The recurrence (\ref{recu1}) yields
\begin{equation}
Q_{n+1}(a) = -a Q_{n}'(a) + (a+n+1)Q_{n}(a).
\label{recu5}
\end{equation}
\noindent
The recursion for the coefficients $b_{n,k}$ follows directly from here.
Morover, it is clear that $b_{n,n}=1$ and $b_{n,0} = n!$. A little 
experimentation suggets that $b_{n,k} = n!/k!$, and this can be established 
from (\ref{recu3}). 
\end{proof}

This proposition amounts to the evaluation of ${\mathbf{3.351.1}}$ in 
\cite{gr}:
\begin{equation}
\int_{0}^{u} x^{n} e^{ax} \, dx = 
\frac{(-1)^{n+1} \, n!}{a^{n+1}} + 
\frac{e^{au}}{a^{n+1}} \sum_{k=0}^{n} \frac{n!}{k!} (-1)^{n-k} u^{k}a^{k}. 
\end{equation}
\noindent
The reader will find a proof of this formula in \cite{moll-gr7}.

\section{A new family of integrals} \label{new} 
\setcounter{equation}{0}

In this section we consider the family of integrals 
\begin{equation}
q_{n}  := \int_{0}^{1} x^{n} e^{-x} \ln x \, dx, \label{qn-def}
\end{equation}
\noindent
and its companion 
\begin{equation}
p_{n}  := \int_{0}^{1} x^{n} e^{-x} \, dx. \label{p-def}
\end{equation}

\begin{lemma}
The integrals $p_{n}, \, q_{n}$ satisfy the recursion
\begin{eqnarray}
p_{n+1} & = & (n+1)p_{n} - e^{-1} \label{recu8} \\
q_{n+1} & = & (n+1)q_{n} + p_{n} \label{recu8a}
\end{eqnarray}
\end{lemma}
\begin{proof}
Integrate by parts. 
\end{proof}

The initial conditions are 
\begin{equation}
p_{0} = 1 - e^{-1} \text{ and } q_{0} = \int_{0}^{1} e^{-x} \ln x \, dx  = \gamma - \text{Ei}(-1). \label{initq}
\end{equation}
\noindent
Here $\gamma$ is Euler's constant defined by
\begin{equation}
\gamma := \lim\limits_{n \to \infty} \sum_{k=1}^{n} \frac{1}{k} - \ln n 
\end{equation}
\noindent
with integral representation 
\begin{equation}
\gamma = \int_{0}^{\infty} e^{-x} \ln x \, dx 
\end{equation}
\noindent
given as ${\mathbf{4.331.1}}$. The 
reader will find in \cite{irrbook}  a proof of 
this identity. The second term in (\ref{initq}) is converted into
\begin{equation}
\int_{1}^{\infty} e^{-x} \ln x \, dx = \int_{1}^{\infty} \frac{e^{-x}}{x} \, dx
\end{equation}
\noindent
and this last form is identified as $\text{Ei}(-1)$, where $\text{Ei}$ is the 
exponential integral defined by
\begin{equation}
\text{Ei}(z) = - \int_{-z}^{\infty} \frac{e^{-x}}{x} \, dx. 
\end{equation}
\noindent
In the current context, the value of $\text{Ei}(-1)$ will be simply one of 
the terms in the initial condition $q_{0}$.

We determine first an explicit expression for $p_{n}$. The recursion 
(\ref{recu8}) shows the existence of integers $a_{n}, \, b_{n}$ such that
\begin{equation}
p_{n} = a_{n} + b_{n}e^{-1},
\end{equation}
\noindent
with $a_{0} =1, \, b_{0} = -1$. From (\ref{recu8}) we obtain
\begin{equation}
a_{n+1} + b_{n+1}e^{-1} = (n+1)a_{n} + \left[ (n+1)b_{n} - 1 \right] e^{-1}.
\end{equation}
\noindent
The irrationality of $e$ produce the system 
\begin{eqnarray}
a_{n+1} & = & (n+1)a_{n}, \text{ with } a_{0}=1, \label{system1} \\
b_{n+1} & = & (n+1)b_{n}-1, \text{ with } b_{0}=-1. \label{system2}
\end{eqnarray}
\noindent
The expression $a_{n} = n!$ follows directly from (\ref{system1}). To solve 
(\ref{system2}), define $B_{n} := b_{n}/n!$ and observe that
\begin{equation}
B_{n+1} = B_{n} - \frac{1}{(n+1)!}, 
\end{equation}
\noindent
that telescopes to 
\begin{equation}
b_{n} = -n! \sum_{k=0}^{n} \frac{1}{k!}. 
\end{equation}

\noindent
We have shown:

\begin{proposition}
The integral $p_{n}$ in (\ref{p-def}) is given by
\begin{equation}
p_{n} = \int_{0}^{1} x^{n} e^{-x} \, dx = 
\frac{n!}{e}  \left( e - \sum_{k=0}^{n} \frac{1}{k!} \right).
\end{equation}
\end{proposition}

We now determine a similar closed-form for $q_{n}$.  The recursion 
(\ref{recu8a}) shows the existence of integers $c_{n}, \, d_{n}, \, f_{n}$ 
such that
\begin{equation}
q_{n} = c_{n} + d_{n}e^{-1} + f_{n} q_{0}. 
\end{equation}
\noindent
In order to produce a system similar to (\ref{system1},\ref{system2}) we will 
assume that the constants $1, \, e^{-1}$ and $q_{0} = -(\gamma + 
\text{Ei}(-1))$ are linearly independent over $\mathbb{Q}$.  Under this
assumption (\ref{recu8a}) produces
\begin{eqnarray}
c_{n+1} & = & (n+1)c_{n} + n!, \label{system3} \\
d_{n+1} & = & (n+1)c_{n}  - n! \sum_{k=0}^{n} \frac{1}{k!},  \label{system4} \\
f_{n+1} & = & (n+1)f_{n}, \label{system5} 
\end{eqnarray}
\noindent
with the initial conditions $c_{0} = 0, \, d_{0} = 0$ and $f_{0} = 1$. 

The expression $f_{n} = n!$ follows directly from (\ref{system5}).  To solve 
(\ref{system3}) and (\ref{system4}) we employ the following result 
established in \cite{moll-gr5}. 

\begin{lemma}
\label{lemrec}
Let $a_{n}, \, b_{n}$ and $r_{n}$ be sequences with $a_{n}, \, 
b_{n} \neq 0$. Assume that
$z_{n}$ satisfies 
\begin{equation}
a_{n}z_{n} = b_{n} z_{n-1} + r_{n}, \, n \geq  1
\label{recu0}
\end{equation}
\noindent
with initial condition $z_{0}$. Then 
\begin{equation}
\label{valuex}
z_{n} = \frac{b_{1}b_{2} \cdots b_{n}}{a_{1}a_{2} \cdots a_{n}} 
\left( z_{0} + \sum_{k=1}^{n} \frac{a_{1}a_{2} \cdots a_{k-1}}
{b_{1} b_{2} \cdots b_{k}} r_{k} \right). 
\end{equation}
\end{lemma}

We conclude that
\begin{equation}
c_{n} = n! \sum_{k=1}^{n} \frac{1}{k},
\end{equation}
\noindent
and 
\begin{equation}
d_{n} = -n! \sum_{k=1}^{n} \frac{1}{k} \sum_{j=0}^{k-1} \frac{1}{j!}. 
\end{equation}
\noindent
The expression for $c_{n}$ shows that they coincide with the  Stirling numbers
of the first kind: $c_{n} = | s(n+1,2)|$. \\

We have established 

\begin{proposition}
The integral $q_{n}$ in (\ref{qn-def}) is given by
\begin{equation}
q_{n} = \int_{0}^{1} x^{n} e^{-x} \ln x \, dx = 
n! \left[ \frac{1}{e} \sum_{k=1}^{n} \frac{1}{k} 
\left( e - \sum_{j=0}^{k-1} \frac{1}{j!}
\right) + q_{0} \right].
\end{equation}
\end{proposition}

\begin{example}
The expressions for $p_{n}$ and $q_{n}$ provide the evaluation of 
$\mathbf{4.351.1}$ in \cite{gr}
\begin{equation}
\int_{0}^{1} (1-x)e^{-x} \ln x \, dx = \frac{1-e}{e},
\end{equation}
\noindent
by identifying the integral as $q_{0}-q_{1}$. The recurrence (\ref{recu8a}) 
shows that 
\begin{equation}
q_{0}-q_{1} = -p_{0} = e^{-1}-1,
\end{equation}
\noindent
as claimed. 
\end{example}

\begin{example}
The evaluation of
$\mathbf{4.362.1}$ in \cite{gr} 
\begin{equation}
\int_{0}^{1} xe^{x} \ln(1-x) \, dx = 
\int_{0}^{1} (1-t)e^{1-t} \ln t \, dt
\end{equation}
\noindent
is achieved by observing that this integral is $e( q_{0} - q_{1}) 
= 1-e$. 
\end{example}
 
\section{A parametric family} \label{sec-param} 
\setcounter{equation}{0}

In this section we consider the evaluation of 
\begin{eqnarray}
P_{n}(a) & := & \int_{0}^{1} x^{n} e^{-ax} \, dx \label{generalp} \\
Q_{n}(a) & := & \int_{0}^{1} x^{n} e^{-ax} \ln x \, dx. \label{generalq}
\end{eqnarray}
\noindent
The integrals $q_{n}$ considered in Section \ref{new} corresponds to the 
special case: $q_{n} = Q_{n}(1)$. 

We now establish a recursion for $Q_{n}$ by differentiating 
(\ref{generalq}).

\begin{lemma}
The integral $Q_{n}(a)$ satisfies the relation
\begin{equation}
Q_{n+1}(a) = -\frac{d}{da} Q_{n}(a).
\end{equation}
\end{lemma}

To obtain a closed-form expression for $Q_{n}(a)$ we need to determine the 
initial condition 
\begin{equation}
Q_{0}(a) = \int_{0}^{1} e^{-ax} \ln x \, dx. 
\end{equation}
\noindent
This is expressed in terms of the {\em incomplete gamma function}
defined in $\mathbf{8.350.1}$ by
\begin{equation}
\Gamma(a,x) := \int_{x}^{\infty} t^{a-1} e^{-t} \, dt.
\end{equation}
\noindent
Observe that $\Gamma(a,0) = \Gamma(a)$, the usual gamma function. 

\begin{lemma}
The initial condition $Q_{0}(a)$ is given by
\begin{equation}
Q_{0}(a) = \int_{0}^{1} e^{-ax} \ln x \, dx = -\frac{1}{a} 
\left( \gamma + \ln a + \Gamma(0,a) \right).
\end{equation}
\end{lemma}
\begin{proof}
The change of variables $t = ax$ yields
\begin{equation}
Q_{0}(a) = \frac{1}{a} \int_{0}^{a} e^{-t} \ln t \, dt - \frac{\ln a}{a} 
\left( 1 - e^{-a} \right). 
\end{equation}
\noindent
Then 
\begin{equation}
\int_{0}^{a} e^{-t} \ln t \, dt = 
\int_{0}^{\infty} e^{-t} \ln t \, dt  - 
\int_{a}^{\infty} e^{-t} \ln t \, dt. 
\end{equation}
\noindent
The first integral is 
\begin{equation}
\int_{0}^{\infty} e^{-t} \ln t \, dt   = - \gamma,
\end{equation}
\noindent
that simply reflects the fact that $\gamma  = - \Gamma'(1)$. 
Integrating by parts yields
\begin{equation}
\int_{a}^{\infty}  e^{-t} \ln t \, dt = e^{-a} \ln a + \Gamma(0,a). 
\end{equation}
\noindent
The formula for $Q_{0}(a)$ is established.
\end{proof}

We now determine a closed-form expression for $P_{n}(a)$ and $Q_{n}(a)$ 
following the procedure employed in Section \ref{new}. 

\begin{lemma}
The integrals $P_{n}$ and $Q_{n}(a)$ satisfy the recursion
\begin{eqnarray}
P_{n+1}(a) = \frac{1}{a} \left( (n+1)P_{n}(a) - e^{-a} \right) \label{recP} \\
Q_{n+1}(a) = \frac{1}{a} \left( (n+1)Q_{n}(a) + P_{n}(a) \right). \label{recQ}
\end{eqnarray}
\noindent
The initial conditions are given by 
\begin{equation}
P_{0}(a) = \frac{1}{a}( 1 - e^{-a}), \text{ and } 
Q_{0}(a) = - \frac{1}{a}( \gamma + \Gamma(0,a) + \ln a ). 
\end{equation}
\end{lemma}
\begin{proof}
Integrate by parts.
\end{proof}

We conclude that we can write 
\begin{equation}
P_{n}(a) = A_{n}(a) - B_{n}(a) e^{-a},
\end{equation}
\noindent
and 
\begin{equation}
Q_{n}(a) = C_{n}(a) - D_{n}(a)e^{-a} - E_{n}(a) 
( \gamma + \Gamma(0,a) + \ln a ). 
\end{equation}

\begin{lemma}
The recursions (\ref{recP}) and (\ref{recQ}) imply that
\begin{eqnarray}
A_{n+1}(a) & = & \frac{1}{a} (n+1) A_{n}(a), \label{recA1} \\
B_{n+1}(a) & = & \frac{1}{a} \left[ (n+1) B_{n}(a) + 1 \right], \nonumber  \\
C_{n+1}(a) & = & \frac{1}{a} \left[ (n+1) C_{n}(a) + A_{n}(a)  \right], \nonumber  \\
D_{n+1}(a) & = & \frac{1}{a} \left[ (n+1) D_{n}(a) + B_{n}(a)  \right], \nonumber  \\
E_{n+1}(a) & = & \frac{1}{a} (n+1) E_{n}(a) \nonumber 
\end{eqnarray}
\noindent
with initial conditions 
\begin{equation}
A_{0}(a) = B_{0}(a) = E_{0}(a) = \frac{1}{a} \text{ and } 
C_{0}(a) = D_{0}(a) = 0. 
\end{equation}
\end{lemma}

These recursion can now be solved as in Section \ref{new} to produce a 
closed-form expression for the integrals $P_{n}(a)$ and $Q_{n}(a)$. We 
employ the notation 
\begin{equation}
H_{n} = \sum_{k=1}^{n} \frac{1}{k}
\end{equation}
\noindent
for the harmonic numbers and 
\begin{equation}
\text{Exp}_{n}(x) = \sum_{k=0}^{n} \frac{x^{k}}{k!}
\end{equation}
\noindent
for the partial sums of the exponential function. 

\begin{theorem}
\label{mainPQ}
Let $a \in \mathbb{R}$ and $n \in \mathbb{N}$. Then 
\begin{equation}
P_{n}(a) := \int_{0}^{1}x^{n}e^{-ax} \, dx = 
\frac{n!}{a^{n+1}} \left[ 1 - e^{-a} \text{Exp}_{n}(a) \right],
\label{formulaP}
\end{equation}
\noindent
and
\begin{equation}
Q_{n}(a) := \int_{0}^{1}x^{n}e^{-ax} \ln x \, dx = 
\frac{n!}{a^{n+1}} \left[  H_{n} - G(a) - e^{-a} \sum_{k=1}^{n} 
\frac{1}{k} \text{Exp}_{k-1}(a) \right],
\label{formulaQ}
\nonumber 
\end{equation}
\noindent
where  $G(a) = -aQ_{0}(a) = \gamma + \Gamma(0,a) + \ln a$. 
\end{theorem}

These expressions provide the evaluations of two integrals in \cite{gr}. 

\begin{example}
Formula ${\mathbf{4.351.2}}$ states that
\begin{equation}
\int_{0}^{1} e^{-ax} (-ax^{2} + 2x) \ln x \, dx = 
\frac{1}{a^{2}} \left[ -1 + (1+a)e^{-a} \right]. 
\end{equation}
\noindent
In order to verify this, observe that the stated integral is 
\begin{equation}
-a \int_{0}^{1} x^{2} e^{-ax}  \ln x \, dx + 
2 \int_{0}^{1} x e^{-ax}  \ln x \, dx  = -a Q_{2}(a) + 2 Q_{1}(a). 
\end{equation}
\noindent
The expressions in Theorem \ref{mainPQ} now complete the evaluation.  
\end{example}

\begin{example}
Formula ${\mathbf{4.353.3}}$ in \cite{gr} gives the value of 
\begin{equation}
I_{n}(a) := \int_{0}^{1} (-ax+ n+1) x^{n} e^{-ax} \ln x \, dx. 
\end{equation}
\noindent
Observe that 
\begin{equation}
I_{n}(a) = -aQ_{n+1}(a) + (n+1)Q_{n}(a),
\end{equation}
\noindent
and using the recursion (\ref{recQ}) we conclude that $I_{n}(a) = 
-P_{n}(a)$. The expression in Theorem \ref{mainPQ} is precisely what 
appears in \cite{gr}. 
\end{example}

We conclude with the evaluation of a series shown to us by 
Tewodros Amdeberhan. Expand the exponential 
term in (\ref{formulaQ}) and integrate term by term to obtain
\begin{equation}
\sum_{k=0}^{\infty} \frac{(-a)^{k}}{k! \, (n+1+k)^{2}} = 
\frac{n!}{a^{n+1}} \left( -\psi(n+1) + \ln a + \Gamma(0,a) +
e^{-a} \sum_{k=0}^{n} \frac{1}{k} \text{Exp}_{k-1}(a) \right). 
\label{series1}
\end{equation}
\noindent
Here 
\begin{equation}
\psi(x) = \frac{\Gamma'(x)}{\Gamma(x)}
\end{equation}
\noindent
is the {\em digamma} function defined in $\mathbf{8.360.1}$ of \cite{gr}. 
the identity 
\begin{equation}
\psi(n+1) = H_{n} - \gamma,
\end{equation}
\noindent
that is a direct consequence of the functional equation $\Gamma(x+1) = 
x \Gamma(x)$ and $\Gamma'(1)= - \gamma$, was used to transform 
(\ref{series1}). 

\medskip

The identity (\ref{series1}) can be used to provide multiple expressions 
for the incomplete gamma function, such as 
\begin{equation}
\int_{a}^{\infty} \frac{e^{-x}}{x} \, dx = 
\sum_{k=0}^{\infty} \frac{(-1)^{k} a^{n+1+k}}{n! \, k! \, (n+1+k)^{2}} + 
\psi(n+1) - \ln a - e^{-a} \sum_{k=1}^{n} \frac{\text{Exp}_{k-1}(a)}{k}, 
\end{equation}
\noindent
and the special case for $n=0$:
\begin{equation}
\int_{a}^{\infty} \frac{e^{-x}}{x} \, dx = 
-\gamma - \ln a + \sum_{k=0}^{\infty} \frac{(-1)^{k} a^{k+1}}
{  (k+1)! \, (k+1)}.   
\end{equation}

\medskip

These issues will be explored in a future publication.

\end{document}